\documentclass[11pt]{article} 

\usepackage[utf8]{inputenc} 

\usepackage{geometry} 
\usepackage[parfill]{parskip} 

\usepackage{graphicx} 
\usepackage{array} 

\usepackage{hyperref,doi}
\usepackage{amsmath,amsthm,amssymb,hyperref,mathtools}
\usepackage{hyperref}
\usepackage{cases}
\usepackage{xcolor}
\usepackage{titlesec}
\usepackage{multirow}
\usepackage{qtree}

\usepackage{caption}
\usepackage{subcaption}

\titleformat{\subparagraph}
    {\normalfont\normalsize\itshape}{\thesubparagraph}{1em}{}
\titlespacing*{\subparagraph}{\parindent}{3.25ex plus 1ex minus .2ex}{.75ex plus .1ex}

\renewcommand{\d}{\partial}

\title{A one-dimensional morphoelastic model for burn injuries:\\\large sensitivity analysis and a feasibility study}
\author{Ginger Egberts $\cdot$ Fred Vermolen $\cdot$ Paul van Zuijlen}

\begin{document}
\maketitle

\begin{abstract}
We consider a one-dimensional morphoelastic model describing post-burn scar contraction. This model describes the displacement of the dermal layer of the skin and the development of the effective Eulerian strain in the tissue. Besides these components, the model also contains components that play a major role in skin repair after trauma. These components are signaling molecules, fibroblasts, myofibroblasts, and collagen.
We perform a sensitivity analysis for many parameters of the model and use the results for a feasibility study. In this study, we test whether the model is suitable for predicting the extent of contraction in different age groups. To this end, we conduct an extensive literature review to find parameter values. 
From the sensitivity analysis, we conclude that the most sensitive parameters are the equilibrium collagen concentration in the dermal layer, the apoptosis rate of fibroblasts and myofibroblasts, and the secretion rate of signaling molecules. Further, although we can use the model to simulate distinct contraction densities in different age groups, our results differ from what is seen in the clinic. 
\end{abstract}

\section{Introduction}
\label{sec:1}
Burns are debilitating, life threatening, and difficult to assess and manage \cite{Lang2019}. Complications after a burn can include, among others, a shock, infection, and long-lasting distress. Further, almost all deeper burns will lead to scarring. The post-burn scars may be immature / mature, atrophic / hypertrophic / keloid, stable / unstable, depigmented (vitiligo) / hyperpigmented, and may turn malignant as well \cite{Goel2010}. In addition, scars are subject to change. For example, an immature scar can mature, and an atrophic scar can become hypertrophic. Post-burn scars are dry and itchy, and need to be prevented from exposure to sunlight. 

One of the common complications in post-burn scars is contraction. Contraction is an active biological process that decreases an area of skin loss in an open wound because of a concentric reduction in the wound's size \cite{Goel2010}. Starting in the proliferative phase of wound healing, wound contraction processes until full scar maturation, after which wound contraction can become active again. Depending on the extent of contraction and the wound dimensions, this can cause limited range-of-motion of joints if the scar is on or in a joint (contracture). This can lead to immobility and is an important indication for scar revision \cite{Egberts2020A}. Scar reconstruction can be necessary.

Burns are so different from other types of wounds that there is a very separate discipline for this class in medical care. There are many different classifications for burns, for example, a burn can be thermal, electrical, or chemical. Burns come with a generalized increase in capillary permeability due to heat effect and damage, and this increase in capillary permeability is not seen in any other type of wound \cite{Tiwari2012}. Burn wound healing consists of three overlapping phases: inflammation (reactive), proliferation (reparative), and maturation (remodeling). During inflammation, the wound is cleaned and prepared for further protection from bacterial infection. Subprocesses in proliferation are reepithelialization, angiogenesis, fibroplasia, and wound contraction. The last phase, in which the scar maturates and attains a balanced structure, can take years. This results in a scar that, on average, has 50\% strength of unwounded skin (within three months), and 80\% on the long-term \cite{Enoch,Young}.

Within the proliferative phase, fibroplasia encompasses the sub-processes that cause the restoration of the presence of fibroblasts and the production of a new extracellular matrix (ECM) in the injured area \cite{KoppenolThesis}. Fibroblasts can differentiate to myofibroblasts which are responsible for pulling forces in the skin and stimulate, like fibroblasts, both the production of the constituents of the new collagen-rich ECM and the release of matrix metalloproteinases (MMPs). The differentiation of fibroblasts to myofibroblasts is stimulated by transforming growth factor $\beta$ \cite{Desmouliere1993}. The fibroblasts, the myofibroblasts and collagen deposition play an important role in the wound contraction. 

For a long time, mathematical models have been developed that simulate the processes involved in wound healing. These models, \cite{Barocas,Dallon1999,KoppenolThesis,McDougall,Olsen1995,Tranquillo} to name a few, predict the behavior of experimental and clinical wounds, and gain insight into which elements of the wound healing response might have a substantial influence on the contraction. The majority of the models can be placed into one of three categories: continuum hypothesis-based models, discrete cell-based models and hybrid models \cite{KoppenolThesis}. One of the subcategories of the continuum hypothesis-based models are the mech\-a\-no-(bio)chemical models. These models together with hybrid models served as a basis for the morphoelastic model that we use in this study and that is developed by Koppenol and Vermolen \cite{Koppenol2017a}.

This model is based on the following principle \cite{Hall}:
the total deformation is decomposed into a deformation as a result of growth or shrinkage and a deformation as a result of mechanical forces. In a mathematical context, one considers the following three coordinate systems: ${\bf X}$, ${\bf X}_e(t)$, and ${\bf x}(t)$, which, respectively, represent the initial coordinate system, the equilibrium at time $t$ that results due to growth or shrinkage, and the current coordinate system that results due to growth or shrinkage and mechanical deformation. Assuming sufficient regularity, the deformation gradient tensor is written by
\begin{equation}
{\bf F} = \frac{\partial {\bf x}}{\partial {\bf X}} = \frac{\partial {\bf x}}{\partial {\bf X}_e} ~ \frac{\partial {\bf X}_e}{\partial {\bf X}} 
= {\bf A} {\bf Z},
\label{frames}
\end{equation}
in which the tensor $Z$ represents the deformation gradient tensor due to growth or shrinkage, and $A$ represents the deformation gradient due to mechanical forces \cite{Goriely2006,Hall,Rodriguez1994}.

The beauty of this model is that we can simulate a permanent deformation, which results from the contraction process in burn wound healing. The main variable in this model is the displacement of the skin ($u$), i.e., the variable that makes us able to determine the surface of the wound, and in later stages, the scar. Besides these results, we can also determine the degree of ‘discomfort’ that the patient experiences. We do this by integrating over the entire tissue, including the undamaged part, yielding the total strain energy density. With this we show to what extent there is an elongation in the entire tissue compared to the situation in which there would be no burn. This elongation, which in principle is simply determined by variations in local displacement, may signal nerves, which may cause the patient to experience a nagging sensation. This leads to discomfort in the patient.

The morphoelastic model compromises many parameters. We know some parameter values, while others are unknown and which we need to estimate. Although Koppenol has provided a great overview of parameter values, parameters vary between patients and even along a piece of skin sample. Hence we are interested in both the sensitivity of the parameters and the feasibility of the model. For the feasibility study, we might choose patient-specific information as input, such as the genetic background, gender, age, the location of the wound on the body, the depth of the wound, or any other. In our search for parameter values, we have seen that a great subset of parameters values is estimated for human skin tissue of different age. For example, the average fibroblast doubling time decreases with age \cite{Simpson2009}, indicating a decrease in fibroblast apoptosis rate with age. We therefore choose to use the patient's age as input for the feasibility study. We summarize the change in parameter values that come with aging, and we use the sensitivity analysis results to vary the parameter values along the domain of computation. In our feasibility study, we define distinct classes of patients of different age for which we simulate many burns. The results show the variations in the relative surface area (RSA) density and the total strain energy density, both for patients of different age. 

The organization of this paper is as follows. Section 2 presents the mathematical model and Section 3 presents the numerical method that we use to approximate the solutions. Subsequently, Section 4 presents the parameter values, Section 5 presents the sensitivity analysis, and Section 6 presents the feasibility study. Finally, Section 7 presents the conclusion and discussion.

\section{The mathematical model}
\label{sec:2}
For the sake of completeness, we present the model that was also used in some of our earlier studies. We model the contraction in burn wounds and scars using partial differential equations that consider the displacement of the dermal layer $({u})$, the displacement velocity of the dermal layer (${v}$), the effective Eulerian strain present in the dermal layer ($\varepsilon$) and the changes in distributions of cells and densities of molecules. This morphoelastic continuum hypothesis-based modeling framework was introduced by Koppenol \cite{Koppenol2017a}. The model follows the assumption from Ramtani \cite{Ramtani1,Ramtani2}, which states that the Young’s modulus of skin depends on the collagen concentration. The model incorporates the evolution in the distributions of fibroblasts ($N$) and myofibroblasts ($M$), the concentrations of \emph{signaling molecules} ($c$) such as cytokines, chemokines, and growth factors, and the collagen concentration ($\rho$). 

We discuss the conservation laws for mass and linear momentum, and the evolution equation that describes how the infinitesimal effective strain changes. Since the domain deforms because of forces exerted by cells, the points in the domain are subject to displacement. We incorporate this local displacement rate by passive convection, which is reflected in the second term in the left-hand side of the equations. 

The signaling molecules play an important role in the immune and inflammation response after wound healing. While these molecules migrate to the wound, they induce directed chemotactic migration of cells. We assume that the diffusion is according to normal Fickian diffusion. Furthermore, we take into account the enhanced secretion by fibroblasts and we assume that a portion of myofibroblasts \cite{Barrientos2008} secrete and consume signaling molecules as well. It is known that MMPs remove signaling molecules from the dermal layer \cite{Mast1996,Sternlicht2001}. The MMP concentration depends on the distribution of (myo)fibroblasts and the collagen concentration \cite{Lindner2012}, and the signaling molecule concentration \cite{Overall1991}. Hence taking the proteolytic breakdown by MMPs into account, we have:
\begin{equation}
\frac{\partial c}{\partial t} + \frac{\partial (c v)}{\partial x} = D_c \frac{\partial^2 c}{\partial x^2} + k_c \left[ \frac{c}{a_c^{I} + c} \right][N + \eta^I M] 
- \delta_c \frac{[N + \eta^{II}M]\rho}{1 + a_c^{II}c} c. \label{PDE_c}
\end{equation}
Here $D_c$ is the Fickian diffusion coefficient of the signaling molecules, $k_c$ is the maximum net secretion rate of the signaling molecules, $\eta^I$ is the ratio of myofibroblasts to fibroblasts in the maximum secretion rate of the signaling molecules, $a_c^{II}$ is the concentration of the signaling molecules that causes the half-maximum net secretion rate of the signaling molecules, $\delta_c$ is the proteolytic breakdown rate parameter of the signaling molecules, $\eta^{II}$ is the ratio of myofibroblasts to fibroblasts in the secretion rate of the MMPs and $1/[1+a_c^{III}c]$ represents the inhibition of the secretion of the MMPs. Note that the MMP balance is assumed to be instantaneous. 
One can find the derivation of this equation in the work of Koppenol and Vermolen \cite{Koppenol2017a}, and partly in the work of Olsen et al. \cite{Olsen1995}. 

We assume that the migration of the cells through the wound bed and surrounding undamaged tissue is determined by random walk and by chemical stimulation by the signaling molecules \cite{Boon2016,Dallon2001,Postlethwaite1987}. This is modeled by a minimal model for chemotaxis \cite{Hillen2008}, and a cell density-dependent Fickian diffusion. We incorporate the proliferation of the cells into the model using logistic growth models, dependent on the signaling molecules (as an activator-inhibitor \cite{Murray2011}), and inhibition because of crowding \cite{VandeBerg1989}. Further, the equations represent differentiation of fibroblasts to myofibroblasts \cite{Tomasek2002}, enhanced by signaling molecules, and apoptosis of the cells:
\begin{multline}
\frac{\partial N}{\partial t} + \frac{\partial (N v)}{\partial x} = -\frac{\partial}{\partial x}\left[- D_F (N+M) \frac{\partial N}{\partial x} + \chi_F N \frac{\partial c}{\partial x}\right]\\
+r_F \left[ 1+\frac{r_F^{\text{max}}c}{a_c^{III}+c} \right][1-\kappa_F (N+M)] N^{1+q} 
- k_F cN - \delta_N N, \label{PDE_N}
\end{multline}
\begin{multline}
\frac{\partial M}{\partial t} + \frac{\partial (M v)}{\partial x} = -\frac{\partial}{\partial x}\left[- D_F (N+M)\frac{\partial M}{\partial x} + \chi_F M \frac{\partial c}{\partial x}\right]\\
+ r_F \left[ \frac{[1+r_F^{\text{max}}]c}{a_c^{III}+c} \right][1-\kappa_F (N+M)] M^{1+q} 
+ k_F cN  - \delta_M M. \label{PDE_M}
\end{multline}
Here $D_F$ represents (myo)fibroblast random diffusion and $\chi_F$ is the chemotactic parameter that depends on both the binding and unbinding rate of the signaling molecules with its receptor, and the concentration of this receptor on the cell surface of the (myo)fibroblasts, $r_F$ is the cell division rate, $r_F^\text{max}$ is the maximum factor of cell division rate enhancement because of the presence of the signaling molecules, $a_c^I$ is the concentration of the signaling molecules that cause half-maximum enhancement of the cell division rate, $\kappa_F(N+M)$ represents the reduction in the cell division rate because of crowding, $q$ is a fixed constant, $k_F$ is the signaling molecule-dependent cell differentiation rate of fibroblasts into myofibroblasts, $\delta_N$ is the apoptosis rate of fibroblasts and $\delta_M$ is the apoptosis rate of myofibroblasts. 
Myofibroblasts only proliferate in the presence of the signaling molecules, hence the difference between equation \eqref{PDE_N} and equation \eqref{PDE_M}. Although myofibroblasts are able to differentiate back to fibroblasts under the influence of Prostaglandin E$_2$ (PGE$_2$) \cite{Garrison2013}, we do not take into account the re-differentiation of myofibroblasts to fibroblasts. 

For collagen, we assume that there is no active transport present, because collagen molecules are large, which reduces their diffusivity. Since collagen is extracellular, it is, next to diffusion, not subject to further active migration mechanisms. Collagen is produced by (myo)fibroblasts \cite{Baum2006}, and enhanced by the secretion by signaling molecules \cite{Ivanoff2005}. The proteolytic breakdown of collagen is like for the signaling molecules related to the MMP concentration:
\begin{equation}
\frac{\partial \rho}{\partial t} + \frac{\partial (\rho v)}{\partial x} = k_\rho \left[ 1 + \left[ \frac{k_\rho^{\text{max}}c}{a_c^{IV} + c} \right] \right] [N + \eta^I M]
-\delta_\rho \frac{[N + \eta^{II}M]\rho}{1 + a_c^{II}c} \rho. \label{PDE_R}
\end{equation}
Here $k_\rho$ is the collagen secretion rate, $k_\rho^{\text{max}}$ is the maximum factor of secretion rate enhancement because of the presence of the signaling molecules, $a_c^{IV}$ is the concentration of the signaling molecules that cause the half-maximum enhancement of the secretion rate of collagen and $\delta_\rho$ is the degradation rate of collagen. A generic MMP affects the reaction kinetics of the signaling molecules and collagen, and is assumed always to be at a local equilibrium concentration. Reasoning for this modeling choice has been to avoid even more complexity and additional unknown parameter values. 

In the equation for the displacement velocity, the Cauchy stress tensor is related to the effective Eulerian strain and displacement velocity gradients by a visco-elastic constitutive relation. 
The body force is generated by an isotropic stress and a pulling force on the extracellular matrix (ECM) by myofibroblasts, proportional to the product of the cell density of the myofibroblasts and a function of the concentration of collagen:
\begin{equation}
\rho_t \left( \frac{\partial v}{\partial t} + 2v\frac{\partial v}{\partial x} \right) = \frac{\partial}{\partial x}\left( \mu\frac{\partial v}{\partial x} + E \sqrt{\rho} \epsilon\right)
+ \frac{\partial}{\partial x}\left( \frac{\xi M\rho}{R^2+\rho^2} \right). \label{PDE_v}
\end{equation}
Here $\rho_t$ represents the total mass density of the dermal tissues, $\mu$ is the viscosity, $E\sqrt{\rho}$ represents the Young's modulus (stiffness), $\xi$ is the generated stress per unit cell density and the inverse of the unit collagen concentration, $R$ is a constant.
 
To incorporate plastic deformation in the equation for the effective strain \eqref{PDE_e}, a tensor-based approach is used that is also commonly used in the context of growth of tissues (such as tumors). We assume that the rate of active change of the effective Eulerian strain is proportional to the product of the amount of effective Eulerian strain \cite{Hall}, the local MMP concentration, and the local signaling molecule concentration:
\begin{equation}
\frac{\partial\epsilon}{\partial t} + v\frac{\partial\epsilon}{\partial x} + (\epsilon-1)\frac{\partial v}{\partial x} = -\zeta\frac{[N+\eta^{II}M]c}{1+a_c^{II}c}\epsilon.  \label{PDE_e}
\end{equation}
Here $\zeta$ is the rate of morphoelastic change (i.e., the rate at which the effective strain changes actively over time).

\subsection{The domain of computation, initial and boundary conditions}
We define the domain of computation by $\Omega_{x,t}=(-L,L)$ with $\overline{\Omega}_{x,t}=[-L,L]$, the closed interval. Similar, we define the wounded area by the subspace $\Omega_{x,t}^w=(-L^w,L^w)$, $L^w<L$. Furthermore, we define the steepness of the boundary of the wound by $s$, that counts for the slope of the constituents on the boundary of the wound. The dimension $x$ is in centimeters and $t$ in days.

\begin{figure}[h]
\centering
\includegraphics[width=\linewidth]{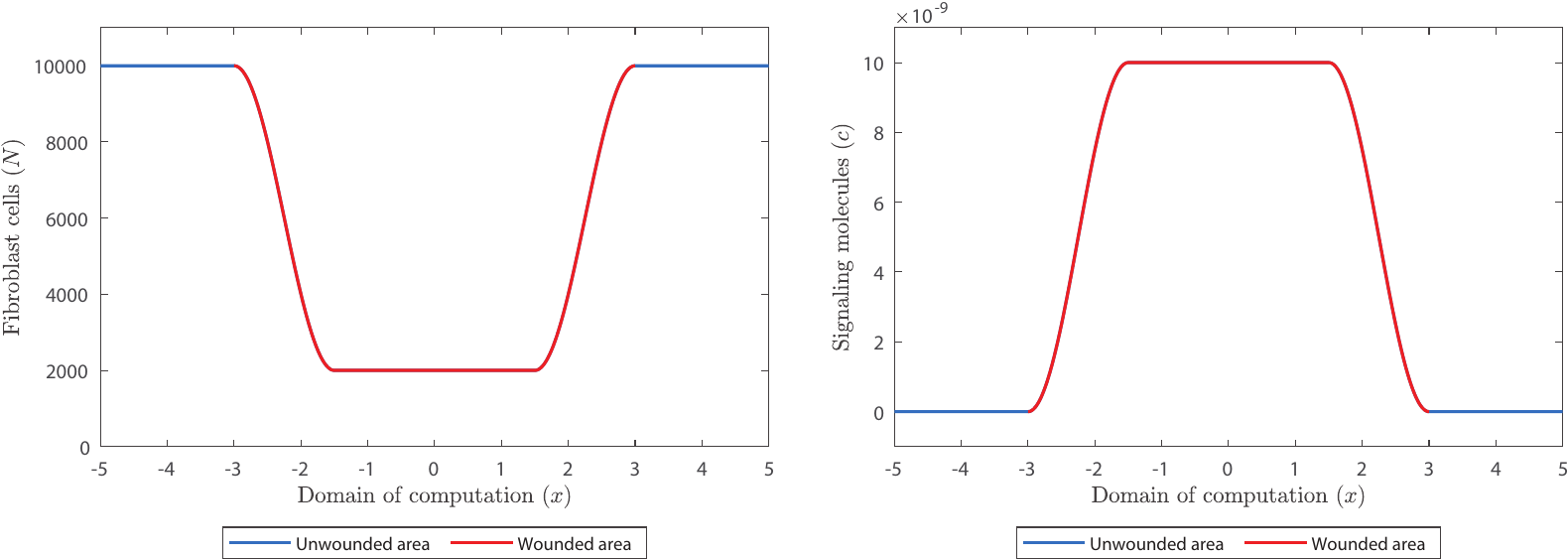}  
\caption{Example of the initial fibroblast density and initial signaling molecule density, with values of the parameters: $L = 5$, $L^w = 3$, $s = 1.5$, $\overline{N}=10^4$, $\tilde{N}=2\times10^3$, $\overline{c} = 0$ and $\tilde{c} = 10^{-8}$.}
\label{fig:InitialConditions}
\end{figure}

We use the following functions for the initial fibroblast density and the initial signaling molecule concentration:
\begin{equation*}
\begin{split}
N(x, 0)=
\left\{\begin{array}{ll}
\overline{N} & \text { if } (*), \\
\frac{\overline{N}+\tilde{N}}{2}+\frac{\overline{N}-\tilde{N}}{2} \sin \left(\frac{\pi}{s}\left(x+\frac{1}{2} s\right)\right) & \text { if } (**), \\
\tilde{N} & \text { if }(***),
\end{array}\right.\\
c(x, 0)=\left\{\begin{array}{ll}
\overline{c} & \text { if } (*), \\
\frac{\overline{c}+\tilde{c}}{2}+\frac{\overline{c}-\tilde{c}}{2} \sin \left(\frac{\pi}{s}\left(x+L^w-\frac{s}{2}\right)\right) & \text { if } (**), \\
\tilde{c} & \text { if }(***).
\end{array}\right.
\end{split}
\end{equation*}
$(*):L^w\leq x \leq-L^w$, $(**):\{L^w + x\leq s, L^w - s\leq x\leq L^w\}$, $(***):-L^w+s \leq x \leq L^w-s$.
Here $\overline{N},\overline{c},\tilde{N},\tilde{c}$ are the fibroblast density and signaling molecule concentration in healthy dermal tissue and in the wound, respectively. Examples of possible initial densities are shown in Figure \ref{fig:InitialConditions}. We use these functions in order to avoid steep changes in the densities. We have assumed that some fibroblasts are present in the wound and used that signaling molecules are present in the wound due to the secretion by for instance macrophages during inflammation.

In line with Koppenol, we assume that initially no myofibroblasts are present, that a fixed collagen concentration is present, and that the displacement of the dermal layer, the displacement velocity and effective Eulerian strain initially are zero. Hence for all $x\in\Omega_{x,0}$:
\begin{equation}\begin{aligned}
M(x, 0)=\overline{M}=0, & \quad \text {and} \quad \rho(x, 0)=\overline{\rho}, \\
u(x, 0)=0, \quad v(x, 0)=0, & \quad \text {and} \quad \varepsilon(x, 0)=0.
\end{aligned}\end{equation}

We impose the following boundary conditions. For all $x\in\d\Omega_{x,t}$ and $t\geq 0$:
\begin{equation}\label{BC}
\begin{aligned}
c(x,t) = 0, &\quad N(x,t)=\overline{N},\quad M(x,t) = 0,\quad v(x,t)=0.
\end{aligned}
\end{equation}
The boundary condition for the displacement velocity follows from the assumption that the boundary of computation is sufficiently far away from the boundary of the wound. We impose no boundary conditions on the collagen concentration, the displacement of the dermal layer and the effective Eulerian strain.

\section{The numerical method}
We approximate the solution to the model equations by the finite-element method using linear basis functions. For more information about this method we refer to Van Kan et al. \cite{VanKan}. We multiply the equations \eqref{PDE_c}-\eqref{PDE_e} by a test function $\varphi(x,t)\in H_0^1$, integrate over the domain of computation $\Omega$ (integration by parts), apply the application of the Gauss' theorem, and apply the Leibniz-Reynold\rq{}s transport theorem. 

To construct the basis functions we subdivide the domain of computation into $n\in\mathbb{N}$ sub-domains $e_p=[x_p,x_{p+1}]$ (i.e., the elements). Let $X_h(t) = \bigcup e_p$ the finite element subspace and $x_j , j \in \{1,\dots,n+1\}$ the vertices of the elements. We choose $\varphi_i(x_j,t)=\delta_{ij}$, $i,j\in\{1,\dots,n+1\}$ as the linear basis functions, where $\delta_{ij}$ denotes the Kronecker delta function. 

Note that the following holds for the chosen subspace $X_h(t) \subset \Omega_{x,t}$: $\frac{\mathrm{D}\varphi_i}{\mathrm{D}t}=0$ for all $\varphi_i$ \cite{Dziuk}. The Galerkin equations are simplified using this property. We solve the Galerkin equations using backward Euler time integration and we use a monolithic approach with inner Picard iterations to account for the non-linearity of the equations. To avoid loss of monotonicity (i.e. oscillations), we use mass lumping. 

We approximate the local displacements of the dermal layer $(u)$ with
\begin{equation}u_{i}^{\tau+1}=u_{i}^{\tau}+\int_{\tau \Delta t}^{(\tau+1) \Delta t} v(s) \mathrm{d} s \simeq u_{i}^{\tau}+\Delta t v_{i}^{\tau+1}.\end{equation}
 
\section{Parameter values}
In our search for parameter values, we used various sources from the literature. One of the most important sources is Koppenol’s thesis, in which several parameter values have been estimated that we did not find in existing literature. In our study we conduct a sensitivity analysis and a feasibility study. In the Appendix one finds the parameter values that we used in this study. In this chapter, we describe how we chose these values.

\paragraph{Equilibrium values}
Taking into account the reaction term for the signaling molecules, and the equilibria $\overline{N},\overline{M}$ and $\overline{\rho}$, the equilibrium density of the signaling molecules should be $\overline{c}=0$ g/cm$^3$.

The estimation of the equilibrium distribution of fibroblasts differs per study. One estimates the number to be about $\mathcal{O}(10^4)$ cells/cm$^3$ \cite{Olsen1995}, and the other estimates the number to be about $\mathcal{O}(10^6)$ cells/cm$^3$ \cite{Miller2003}. The estimation of the number of cells also differs for the papillary and the reticular dermis, where there exist much more fibroblasts in the papillary dermis \cite{Harper1979,Randolph1998}. In our simulations we have seen that the model works best with the equilibrium distribution of $\mathcal{O}(10^4)$ cells/cm$^3$. We note that some other parameter values ($\delta_c,\delta_\rho$) depend on the chosen order for $\overline{N}$, since we need to take into account the density of MMPs. Furthermore, research has found that among ages 1-10 the number of fibroblasts is nearly twice as high as in any other postnatal age group \cite{Gunin}. Therefore we choose the mean value $\overline{N}=10^4$ cells/cm$^3$ and let the value decrease with age.

The number of myofibroblasts present in the skin depends on the condition of the skin. Myofibroblasts result from the differentiation of fibroblasts. We assume that myofibroblasts are not present at the beginning of proliferation, since healthy skin contains almost no myofibroblasts. Hence $\overline{M} = 0$ cells/cm$^3$.

Olsen et al. estimate the equilibrium collagen concentration as follows. Roughly 75\% of the 15\% of other substances than water and fat in 1 ml of human dermal tissue is collagen \cite{Olsen1995}. This yields $\overline{\rho} \approx 0.75 \times 0.15$ g ml$^{-1} = 0.1125$ g ml$^{-1}$. Furthermore, in human skin the collagen content decreases at about 2\% per year \cite{Farage2015}. Therefore we choose the mean value $\overline{\rho} = 0.1125$ g/cm$^3$ and let the value decrease exponentially with age.

\paragraph{Initial values}
Due to the supply of growth factors in the inflammatory phase, the initial signaling molecule concentration is unequal to zero. The value should not exceed 15-50 ng ml$^{-1}$ \cite{Olsen1995}, and is therefore chosen to be $\tilde{c}= 10^{-8}$ g/cm$^3$.

The thermal injury causes sudden death of cells. The dead cells lose their solid integrity, which causes the release of cytokines. These cytokines trigger the immune response where several types of immune cells clear up the debris and release signaling molecules, which trigger the fibroblasts to migrate to the damaged region. Since we simulate from the onset of proliferation, we assume that several fibroblasts are present. We let this number be 20 percent of the equilibrium number. So the mean value is $\tilde{N}=2\times10^3$ cells/cm$^3$. 

\paragraph{Flux values}
Sillman et al. vary the migratory rates of fibroblasts depending on the experimental medium used: in serum-containing medium the average velocity was as low as 0.23 mm/min, while in serum-free keratinocyte medium the average velocity was as high as 0.36 mm/min \cite{Sillman2003}. Hence, in serum-containing medium the rate was $7.6176 \times 10^{-7}$ cm$^2$/day and in serum-free keratinocyte medium the rate was $1.86624 \times 10^{-6}$ cm$^2$/day. All the reported values together yield a mean value of $1.3247 \times 10^{-6}$ cm$^2$/day and standard deviation $3.7823 \times 10^{-8}$ cm$^2$/day. However, other estimates are $1.44 \times 10^{-5}$ cm$^2$/day and $1.2 \times 10^{-5}$ cm$^2$/day \cite{Olsen1995} \& \cite{Simpson2017}. We therefore estimate the value of $D_F \approx 10^6$ cm$^5$/(cells day). Furthermore, we assume that the diffusion of (myo) fibroblasts decreases with age.

For the chemotactic parameter, we adopt $\chi_F = 2 \times 10^{-3}$ cm$^5$/(cells day) from Murphy et al. \cite{Murphy2012}. For the diffusion parameter, we adopt $D_c \approx 2.88\times10^{-3}$ cm$^2$/day from Haugh \cite{Haugh2006}. Furthermore, we assume that the diffusion of signaling molecules decreases with age.

\paragraph{Chemical kinetics values}
Olsen et al. relate the inhibitor of TGF-$\beta$ to the initial concentration of the growth factors so that $a_c^{I}=10^{-8}$ g/cm$^3$ \cite{Olsen1995}. We adopt this value.

Myofibroblasts produce roughly twice the collagen that is synthesized by fibroblasts \cite{Rudolph1991}. Hence the constant $\eta^I=2$.

The half-life of TGF-$\beta$ is about 2 minutes \cite{Wakefield1990}, and the half-life of PDGF is about 2 minutes as well \cite{Bowen-Pope1984}. So signaling molecules have a decay rate of $-\log(0.5^{24\times60/2})$\\ $\approx 499$/day. However, Olsen et al. decrease the value for two reasons: not all signaling molecules may bind, for example because of insufficient levels of binding protein present at the wound site, and the bound complex may be recognized by (myo)fibroblasts leasing to internalized and metabolized signaling molecules \cite{Olsen1995}. Therefore, the estimated decay rate is $0.5$/day. Other estimates for TGF-$\beta$ are $0.462-0.693$/day \cite{Javierre2009} and $0.354$/day \cite{Murphy2012,Yang1999}. Given our equilibrium parameter values, the MMP density has order of magnitude $\mathcal{O}(\overline{N})\times\mathcal{O}(\overline{\rho})=\mathcal{O}(10^3)$. Hence taking care of the equilibrium dimensions of the model, we end up with a range of $(3.54-6.93)\times10^{-4}$ cm$^6$/(cells g day). We take the value $\delta_c = 5 \times 10^{-4}$ cm$^6$/(cells g day). 

From our previous stability analysis it follows that $k_c\leq\delta_c\overline{\rho}a_c^{I}$ \cite{Egberts2020A}. Given the parameter values, we set $k_c = 3 \times 10^{-13}$ g/(cells day).

We estimate the constant $\eta^{II}=0.45$, which is a small deviation from the constant estimated in \cite{KoppenolThesis}.

Overall et al. estimate $a_c^{II} = (2-2.5) \times 10^8$ cm$^3$/g \cite{Overall1991}. We choose the lower limit, hence $a_c^{II} = 2 \times 10^8$ cm$^3$/g. Furthermore, the production of MMPs increases with age \cite{Ashcroft1997}. Given the equation for the MMPs \cite{Koppenol2017a}, that is 
\begin{equation}
g(N,M,c,\rho)=\frac{[N+\eta^{II}M]\rho}{1+a_c^{II}c},
\end{equation}
to let the production of MMPs increase, we must decrease the inhibiting factor $a_c^{II}$. Hence, we let the inhibition factor $a_c^{II}$ decrease with age.

Cell doubling time can be calculated using the growth rate (amount of doubling in one unit of time) in the following way: doubling time = ln(2)/growth rate. The average doubling time for fibroblasts is approximately 18-20 h \cite{Alberts1989,Gosh}. This gives the range for the proliferation rate of $0.832\leq r_F\leq 0.924$. We choose the upper limit, hence $r_F = 0.924$ cm$^{3q}$/(cells$^q$ day). 
Furthermore, the percentage of PCNA-positive fibroblasts decreases with age, and PCNA can be considered as a marker for the proliferating cells \cite{Gunin}. We therefore let the cell division decrease with age.

TGF-$\beta$ increases fibroblast proliferation by 2-3 times \cite{Strutz2001}. We choose the upper limit, hence $r_F^{\text{max}}=2$. 

The chemical concentrations required to enhance fibroblast proliferation are somewhat higher than those for chemotactic responses \cite{Grotendorst1992}. Experimental evidence indicates that half-maximal enhancement corresponds to concentrations about 10 ng per ml \cite{Olsen1995}. We adopt this value and take $a_c^{III}=10^{-8}$ g/cm$^3$.

The carrying capacity of fibroblasts is known to be approximately $\kappa_F=10^{-6}$ cm$^{3}$/cells \cite{VandeBerg1989}. This value is adopted in this study. Furthermore, skin becomes thinner with age and therefore we assume that crowding occurs faster in elderly. Hence we let the division rate reduction value increase with age. 

We need to have a stable chemical reaction in case the cell distributions and molecules densities are in equilibrium. The constant $q$ gives us the opportunity to have a stable reaction in equilibrium for equation \eqref{PDE_N}. Given the equilibria, solving for $q$ yields:
\begin{equation}\label{q}
q=\frac{\log(\delta_N)-\log(r_F(1-\kappa_F\overline{N}))}{\log(\overline{N})}.
\end{equation}

In Desmouli\`{e}re et al. \cite{Desmouliere1993}, culturing fibroblasts in the presence of TGF-$\beta$ increased the percentage of cells expressing $\alpha$-SMA from 7.5\% to 45.3\%, representing an activation of 37.8\% of myofibroblast type cells. This experiment occurred over a one week period, with a TGF-$\beta$ dose of 5-10 ng per ml. Suppose the activation of myofibroblasts follows a linear equation. Then given $y(7) = 7a = 0.378$, we have $a=0.054$/day. A dose of 5-10 ng per ml yields $0.054/10 \times 10^{-9}$ and $0.054/5 \times 10^{-9}$ cm$^3$/(g day), giving the range $5.4 \times 10^6\leq k_F \leq 1.08 \times 10^7$ cm$^3$/(g day). We choose the upper limit. 
Furthermore, Simpson et al. demonstrated a failure of fibroblast-myofibroblast differentiation and showed that this is associated with in vitro aging \cite{Simpson2009}. Hence we let the differentiating parameter decrease with age.

The average fibroblast doubling time ($DT$) ranges between 18-20 h \cite{Alberts1989,Gosh}, and the average lifespan of fibroblasts varies between 40 and 70 population doublings ($PD$) \cite{Azzarone1983,Moulin2011}. Using the formula
\begin{equation} \delta_N = (\ln 2)/(PD \times DT/24), \end{equation}
we end up with the range $0.0119\leq \delta_N\leq 0.0231$. We choose the value $\delta_N=0.02$/day and let this value decrease with age, since on average, the doubling time of fibroblasts decreases with age \cite{Simpson2009}. 

The apoptosis rate of myofibroblasts was estimated in a previous study for hypertrophic scars \cite{Koppenol2017b}. Within this study, it was found that a value of $\delta_M=0.002$ /day corresponds to hypertrophic scars and that a value of $\delta_M = 0.06$ /day corresponds to normal scars. Further averages are: 8.85\% for normal scars and 1.06\% for hypertrophic scars \cite{Moulin2003}. Combination of these results yield the range $0.06\leq\delta_M\leq 0.0885$ for normal scars and $0.0106\leq \delta_M\leq 0.02$ for hypertrophic scars. For our study we use the lower value $\delta_M=0.06$/day for normal scars.

The secretion rate of collagen $k_\rho$ gives us the opportunity to have a stable reaction in equilibrium for equation \eqref{PDE_R}. Given the equilibria, solving yields 
$k_\rho=\delta_\rho\overline{\rho}^2$.

The synergistic effects of growth factors may accelerate collagen biosynthesis up to tenfold \cite{Olsen1995}. Hence $k_\rho^{\text{max}}=10$. 

Data suggests that the half-maximal enhancement of collagen synthesis occurs at TGF-$\beta$ concentrations of the order of 1 ng per ml \cite{Roberts1986}. We adopt this value, hence $a_c^{IV}=10^{-9}$ g/cm$^3$.

For the decay rate of collagen, we let $\delta_\rho = 6\times10^{-6}$ cm$^6$/(cells g day) \cite{Koppenol2017b}. Furthermore, the collagen turnover decreases with age \cite{Farage2015}. Hence we let the proteolytic breakdown of collagen decrease with age.

\paragraph{Mechanical values}
Koppenol et al. estimated the viscosity value to be of order $O(10^2)$ for the two-dimensional morphoelastic model \cite{Koppenol2017a}. In our previous study, the stability analysis showed that $\mu\geq\frac{\sqrt{\rho E}}{\pi}$ must hold for the one-dimensional morphoelastic model \cite{Egberts2020A}. Given other parameter values we can adopt the value $\mu=100$ (N day)/cm$^2$. Furthermore, since the viscosity is constant for patients up to their 40s and increases a little after turning 40 \cite{Xu}, we let the viscosity increase with age. 

We estimate that the constant $E$ in the Young's Modulus $E\sqrt{\rho}$ is $350$ N/((g cm)$^\frac12$) for the one-dimensional morphoelastic model and let this value increase with age \cite{Pawlaczyk,Pond}.

For the parameters in the body force, we adopt $\xi=4.4\times10^{-2}$ (N g)/(cells cm$^2$) \cite{Maskarinec} \& \cite{Wrobel}, $R=0.995$ g/cm$^3$, and $\zeta=(0-9)\times10^2$ cm$^6$/(cells g day). We set $\zeta=4 \times 10^2$/(cells g day) and let this value increase with age, because the skin's ability to recover after stretching decreases over lifetime \cite{Krueger}. 

Last, but not least, $\rho_t=1.09$ g/cm$^3$ for human skin \cite{Wrobel2009}. We asumme this density does not change with age.

\section{Sensitivity analysis}
\label{sec:4}
The model contains 34 parameters of which we vary the following independent 30 to study the sensitivity of these parameters:
\begin{itemize}
\item the equilibria $\overline{N}$ and $\overline{\rho}$, and the initial conditions $\tilde{N},\tilde{c}$ and $\tilde{\rho}$;
\item the apoptosis rates $\delta_N$ and $\delta_M$, and the decay rates $\delta_c$ and $\delta_\rho$;
\item the parameters responsible for enhancement of cell division and molecule secretion $a_c^I,a_c^{III}$ and $a_c^{IV}$, and the inhibition of MMP secretion $a_c^{II}$;
\item the ratios from myofibroblasts to fibroblasts $\eta^I$ and $\eta^{II}$, and chemokine dependent differentiation rate $k_F$;
\item the diffusion and chemotaxis rates $D_F,D_c$ and  $\chi_F$;
\item the proliferation and secretion rates $r_F$ and $k_c$, and the maximum factors $r_F^{\text{max}}$ and $k_\rho^{\text{max}}$;
\item the crowding factor $\kappa_F$;
\item the parameters $\xi$ and $R$ that influence the force;
\item the viscosity $\mu$, Young's-Modulus factor $E$, morphoelastic factor $\zeta$, and the total mass density of dermal tissues $\rho_t$.
\end{itemize}
We also vary the length of the initial wound $L^w$. The analysis is organized as follows. For each chosen parameter we vary the value by decreasing or increasing it by $\pm 0,5,10,15,20,25\%$. This means we perform 341 simulations. Namely, for each parameter we perform 11 simulations while leaving the values of the other parameters at the mean value. The mean values are given in Tables \ref{Parm1} - \ref{Parm4}, where in Tables \ref{Parm3} and \ref{Parm4} the mean values are given in the third column ($\mu^2$).

We are interested in the contraction during wound healing and scar formation, and in the total strain energy density. The latter is assumed to be a measure for the discomfort that the patient experiences and is defined by
\begin{equation}\label{strainenergy}
E_\varepsilon(t) = \int_{-L}^{L}\frac12 E\sqrt{\rho(x,t)}\varepsilon(x,t)^2\,\mathrm{d}x 
= \int_0^L E\sqrt{\rho(x,t)}\varepsilon(x,t)^2\,\mathrm{d}x.
\end{equation}
Here we used the symmetry of the domain. 

The results show the \emph{minimum of the relative surface area (RSA$_{min}$)} (i.e., maximum contraction value) in a time period of one year, the \emph{day on which the minimum relative surface area is reached (RSA$_{day}$)} (i.e., the day after which the wound/scar retracts), the \emph{relative surface area on day 365 (RSA$_{365}$)}, the \emph{maximum value of the strain energy density (SED$_{max}$)}, and the \emph{day on which the maximum value of the strain energy density is reached (SED$_{day}$)} (i.e., the day after which the patient experiences a reduction in discomfort due to the internal stress in the skin). In our simulations we used $L=10$ cm for the boundary of the domain of computation and the mean value $L^w=3.6$ cm for the wound.

Given the values in $r\in RSA_{\{min,day,365\}}$ and $r\in SED_{\{max,day\}}$ for a variation $j\in\{\pm25\%\}$, we compute the \emph{z-scores} for the parameter $i\in\{\overline{N},\dots,$ $L^w\}$. The basic $z$-score for a sample is $z = (x-\overline{x})/s_x$, where $\overline{x}$ is the sample mean and $s_x$ is the sample standard deviation.

We define the measure for sensitivity by summing over the absolute values of the $z$-scores as follows:
\begin{equation}\label{sens_measure}
\mathcal{S}_i^r = \sum_{j} \left| z_{ij}^r \right|,
\end{equation}
where $z_{ij}^r$ is the $z$-score of the data in $r$ for parameter $i$ in variation $j$.
For example, $z_{\delta_N,15\%}^{RSA_{365}}$ represents the $z$-score of the relative surface area on day 365 for parameter $\delta_N$ in the simulation where the value for $\delta_N$ is increased with 15\%.

\begin{table}
\centering
\begin{tabular}{l|lllll|l}
Parameter 			& $\mathcal{S}^{RSA_{min}}$ & $\mathcal{S}^{RSA_{day}}$ & $\mathcal{S}^{RSA_{365}}$ & $\mathcal{S}^{SED_{max}}$ & $\mathcal{S}^{SED_{day}}$ & $\mathcal{S}^\text{total}$\\ \hline
$\overline{N}$			& 9.207 & 8.835 & 9.496 & 10.349 & 7.58 & 45\\
$\overline{\rho}$		& 20.351 & 33.916 & 15.383 & 29.19 & 30.21 & 129\\
$\tilde{N}$			& 3.347 & 2.959 & 1.979 & 3.813 & 3.321 & 15\\
$\tilde{c}$			& 3.687 & 3.015 & 2.065 & 2.721 & 4.544 & 16\\ 
$\tilde{\rho}$			& 1.418 & 3.107 & 1.765 & 2.991 & 2.048 & 11\\
$\delta_N$			&22.549 & 11.747 & 23.891 & 17.851 & 15.05 & 91\\
$\delta_M$			& 20.52 & 11.794 & 23.787 & 17.245 & 11.423 & 85\\
$\delta_c$			& 12.001 & 9.09 & 11.558 & 5.1 & 16.287 & 54\\
$\delta_\rho$			& 8.905 & 15.125 & 6.879 & 13.054 & 13.552 & 58\\
$a_c^I$				& 4.496 & 9.211 & 5.731 & 10.169 & 2.257 & 32\\
$a_c^{II}$			& 4.828 & 2.057 & 4.779 & 2.823 & 4.499 & 19\\
$a_c^{III}$			& 8.929 & 10.753 & 8.912 & 11.091 & 8.917 & 49\\
$a_c^{IV}$			& 2.306 & 3.107 & 1.529 & 2.991 & 3.037 & 13\\
$\eta^I$				& 8.175 & 3.364 & 8.227 & 4.821 & 8.248 & 33\\
$\eta^{II}$			& 6.551 & 5.713 & 9.24 & 6.304 & 5.607 & 33\\
$k_F$				& 1.329 & 7.903 & 2.324 & 8.49 & 2.259 & 22\\
$D_F$				& 1.954 & 2.057 & 1.604 & 1.867 & 2.427 & 10\\
$D_c$				& 2.154 & 3.015 & 1.595 & 2.882 & 2.297 & 12\\
$\chi_F$				& 1.817 & 3.107 & 1.598 & 2.991 & 2.512 & 12\\
$r_F$				& 4.649 & 1.486 & 3.69 & 2.301 & 4.07 & 16\\
$k_c$				& 16.993 & 15.121 & 14.234 & 10.024 & 22.49 & 79\\
$r_F^{\text{max}}$	& 17.546 & 9.879 & 18.152 & 13.054 & 12.097 & 71\\
$k_\rho^{\text{max}}$	& 8.268 & 12.692 & 6.455 & 11.591 & 12.036 & 51\\
$\kappa_F$ 			& 1.89 & 3.107 & 1.573 & 2.991 & 2.551 & 12\\
$\xi$				& 14.892 & 2.057 & 14.571 & 1.18 & 11.62 & 44\\
$R$					& 20.174 & 10.753 & 21.508 & 13.568 & 13.648 & 80\\
$\mu$				&1.782 & 3.107 & 1.6 & 2.991 & 2.423 & 12\\
$E$					&13.306 & 4.157 & 6.397 & 4.94 & 8.485 & 37\\
$\zeta$				& 3.845 & 5.098 & 1.728 & 2.991 & 10.851 & 25\\
$\rho_t$				& 1.817 & 3.107 & 1.598 & 2.991 & 2.512 & 12\\
$L^w$				& 2.038 & 9.173 & 11.521 & 9.795 & 4.011 & 37 
\end{tabular}
\caption{Sensitivity of the varied parameters in terms of $z$-scores}
\label{zscores}
\end{table}

Table \ref{zscores} gives an overview of the sensitivity values in terms of $z$-scores for the 31 parameters that we varied. In the last column we rounded the sum of the values. From this table we can see that the parameter that represents the equilibrium collagen concentration ($\overline{\rho}$) with score 129 is the most sensitive. It is therefore interesting to study the equilibrium collagen concentrations in human skin, since collagen concentrations decrease with age \cite{Farage2015} and we use this value for our age study in the next section. Parameters that are the least sensitive are the diffusion rate of (myo)fibroblasts ($D_F$) with score 10 and the initial collagen concentration ($\tilde{\rho}$) with score 11. Concerning the diffusion $D_F$, we must note that the mean value is of order $\mathcal{O}(10^{-6})$, which is different from the order used by Koppenol \cite{KoppenolThesis}, where it is of order $\mathcal{O}(10^{-7})$. This may lead to stating that variations of this parameter have no major impact on the simulations, while in a different geometry it might be much more sensitive. Concerning the initial collagen concentration ($\tilde{\rho}$), we note that the value is varied when the equilibrium collagen concentration is fixed to the mean value ($\overline{\rho}=0.1125$). In case the equilibrium collagen concentration $\overline{\rho}$ is varied, the parameter for the initial collagen concentration is fixed to 20\% of 0.1125, which is the mean value of $\overline{\rho}$ and not the variation.

Other parameters that seem significantly sensitive are the apoptosis of fibroblasts ($\delta_N$) with score 91, the apoptosis of myofibroblasts ($\delta_M$) with score 85, the constant $R$ that influences the force with score 80, and the secretion of signaling molecules ($k_c$) with score 79. However, we must note that Koppenol estimated the mean values for $\delta_M$, $R$ and $k_c$ for a two-dimensional setting in which he used other values for the parameters \cite{KoppenolThesis}. Furthermore, the value of $k_c$ is based on a stability constraint $k_c\leq\delta_c a_c^{I}\overline{\rho}$, that was found in our previous stability analysis \cite{Egberts2020A}, and since the parameter for the equilibrium collagen concentration is sensitive, it is not a surprise that this secretion parameter is also sensitive. The value for the secretion rate of signaling molecules is not that straightforward. The secretion rate of cytokines is different from the secretion rate of growth factors and yet we model these together in one variable $c$ representing signaling molecules. To prevent the model from unnecessary complicated computations, we continue modeling with this simplification and bear in mind the sensitivity of the parameter $k_c$.  It is therefore also interesting to study the rates of apoptosis of fibroblasts, since the doubling time of fibroblasts decreases with age \cite{Simpson2009}, and we use this value for our age study in the next section.

\begin{figure}
\begin{subfigure}{.48\linewidth}
  \centering
  \includegraphics[width=\linewidth]{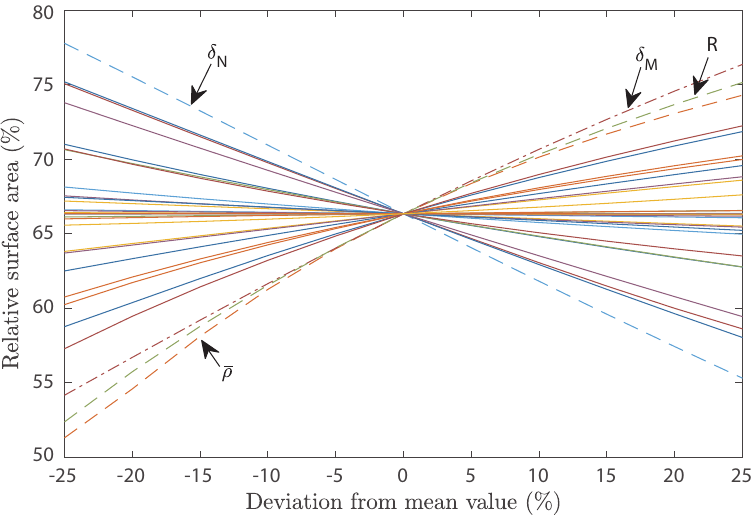}  
  \caption{}
  \label{fig:rsa-min}
\end{subfigure}\hfill
\begin{subfigure}{.48\linewidth}
  \centering
  \includegraphics[width=\linewidth]{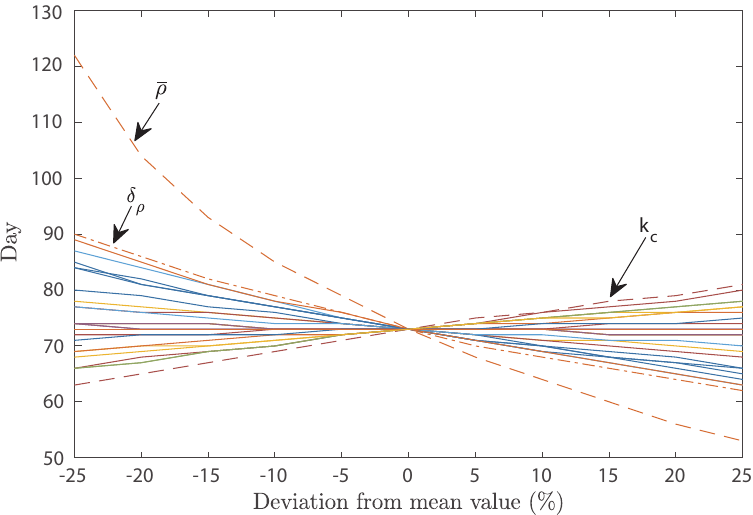}  
  \caption{}
  \label{fig:rsa-day}
\end{subfigure}
\\
\begin{subfigure}{\linewidth}
  \centering
  \includegraphics[width=0.48\linewidth]{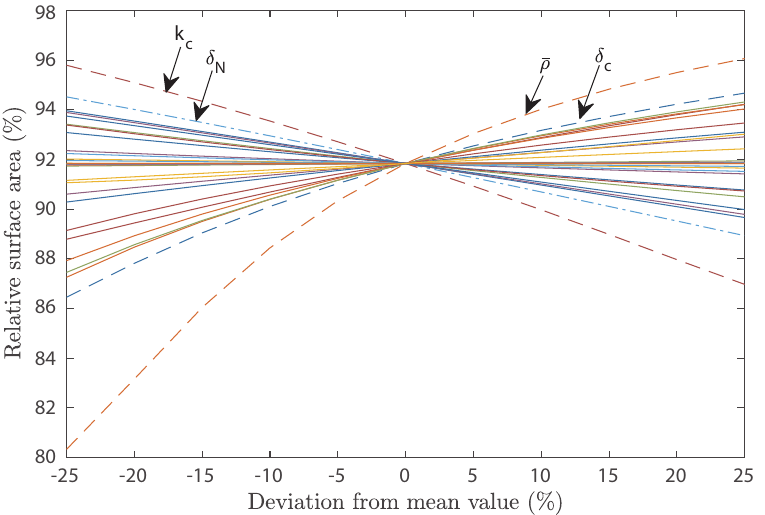}  
  \caption{}
  \label{fig:rsa-365}
\end{subfigure}
\caption{Effects of the variations in parameters for the relative surface area. Shown are the effects on the minimum of the relative surface area (a), the effects on the day on which the minimum of the surface area is reached (b), and the effects on the relative surface area on day 365 (c)}
\label{fig:rsa}
\end{figure}

\begin{figure}
\begin{subfigure}{.48\linewidth}
  \centering
  \includegraphics[width=\linewidth]{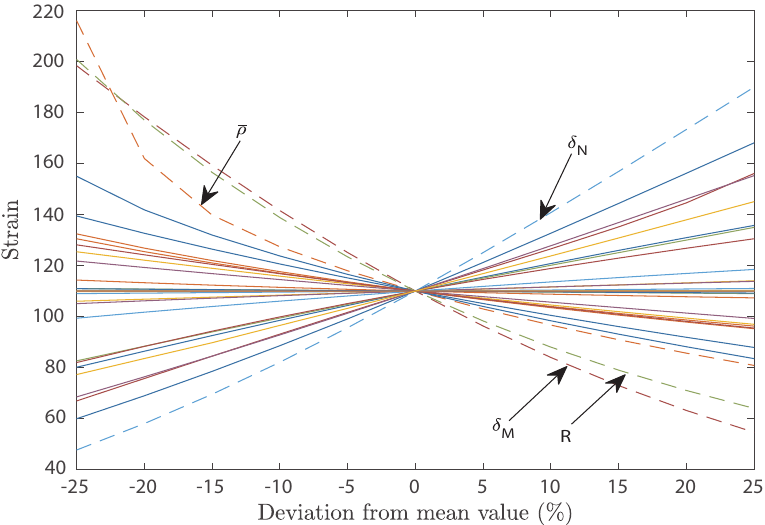}  
  \caption{}
  \label{fig:sed-max}
\end{subfigure}\hfill
\begin{subfigure}{.48\linewidth}
  \centering
  \includegraphics[width=\linewidth]{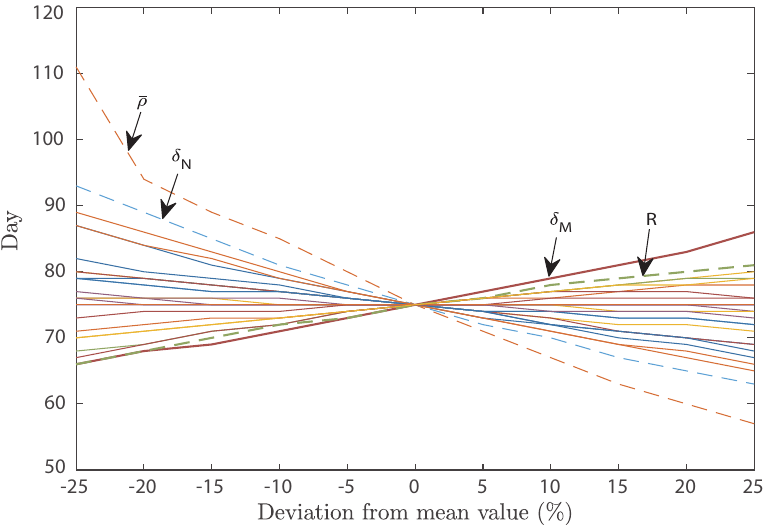}  
  \caption{}
  \label{fig:sed-day}
\end{subfigure}
\caption{Effects of the variations in parameters for the total strain energy density. Shown are the effects on the maximum of the total strain energy density (a) and the effects on the day on which the maximum of the total strain energy density is reached (b)}
\label{fig:sed}
\end{figure}

To get a visual insight into the sensitivity of the parameters, we present the effect of the variations on the parameters on both the relative surface area density of the wound/scar, and the strain energy density in Figures \ref{fig:rsa} \& \ref{fig:sed}. To have a clear distinction between the sensitivity of the parameters, there is no legend, and we labeled the most important lines with different styles. Within the figure, the sub-figures show the effects on the minimum of the relative surface area (a), the effects on the day on which the relative surface area reaches minimum (b), the effects on the relative surface area on day 365 (c), the effects on the maximum of the total strain energy density (d), and the effects on the day on which the total strain energy density reaches maximum (e). 

From Figures \ref{fig:rsa-min} \& \ref{fig:rsa-day}, we see that the parameters that have the most influence on the minimum of the relative surface area and the day on which the relative surface area reaches minimum, are the apoptosis rates of the (myo)fibroblasts, the equilibrium collagen density, and the constant $R$ that together with myofibroblasts and collagen is an important parameter in the body force in equation \eqref{PDE_v}. In particular, a reduction in the fibroblast cell-density yields a reduction in the production of collagen, whereby a reduction in the skin’s stiffness ($E\sqrt{\rho}$). 
Note that a smaller relative surface area corresponds with more intense contraction. Hence, from Figure \ref{fig:rsa-min}, we see that there is less contraction for a smaller fibroblast apoptosis rate, a larger myofibroblasts apoptosis rate, a larger collagen equilibrium density, and a larger value for $R$. The latter three are no surprise, given the body force term in equation \eqref{PDE_v}. Considering the fibroblast apoptosis rate, a reduction in apoptosis rate means more cells survive, hence a relative increase in proliferation. We can see this from equation \eqref{q}: a smaller fibroblast apoptosis rate directly correlates with a smaller value for $q$, and given equation \eqref{PDE_N}, we see that if $q$ becomes smaller, the production of fibroblasts increases. This also relates to the skin’s stiffness, as in an increase in the fibroblast cell-density yields an increase in collagen production, and hence an increase in stiffness. The tissue's strengh therefore increases, and the effective strain decreases, i.e. $\lVert\varepsilon\rVert^2$ becomes smaller.
The result is less contraction during proliferation (Fig. \ref{fig:rsa-min}) and maturation (Fig \ref{fig:rsa-365}), where the effect during proliferation is larger. Further, from Figures \ref{fig:sed-max} \& \ref{fig:sed-day} we conclude that a decrease in fibroblast apoptosis rate results in less discomfort, over a longer period in time. 

In Figures \ref{fig:rsa-day}, the equilibrium collagen density stands out. We see that an increase in the equilibrium collagen density results almost linearly in a significant reduction in the time where the contraction is maximal. In contrast, a decrease of 25\% relative to mean values results in a more exponential-like increase in the maximal contraction time. This property also relates to the total strain energy density, which we see in Figures \ref{fig:sed-max} \& \ref{fig:sed-day}. From the first we see that the discomfort that a patient can experience increases more intensively for a 5 to 20\% reduction in the myofibroblast apoptosis rate and the body force-inhibiting constant $R$. If the equilibrium collagen density decreases with 25\%, then the effect is larger. In reality, it is not likely to change the equilibrium collagen density, however, we can use collagen dressings and, for example, vitamin C supplements to reduce contraction and healing time.

Figures \ref{fig:rsa-day} \& \ref{fig:rsa-365} also feature the secretion rate of the signaling molecules. From the figures we see that a lower signaling molecule secretion rate reduces both the period and intensity of contraction. Especially during maturation, when inflammatory responses are not favorable. 

Taken together, targeting contraction intensity in the proliferative phase of wound healing is most likely effective in case fibroblast survival and collagen density are considered. Targeting contraction, and contractures, during maturation is more likely to be effective when inhibition of signaling molecules and collagen production are considered.

We use the results from the sensitivity analysis to perform a feasibility study in the next chapter.

\section{Feasibility study of modeling age dependent scar/skin contraction}
\label{sec:5}
To study the feasibility of age-dependent uncertainty quantification, we focus on the effect of aging of skin on contraction, the final contracture, the total strain energy density, and the maximum of the total strain energy density. Just like any other organ, aging also affects the skin. Aging has a delaying effect on wound healing and immune responses. Intrinsic aging is the effect of generic and internal influences, such as hormones or metabolic substances. Extrinsic aging is the effect of external influences, such as UV radiation and environmental toxins \cite{Wiegand}. Clear general signs of aging are wrinkles, sagging skin and pigmentary irregularities, and increased tendencies to injuries and the faster opening of healing wounds. These symptoms result from physiological changes such as decreased cell replacement rate. We review various sources from literature to find suitable values for the parameters of the model. In this way we can perform simulations for patients of different ages. Based on the results found, we have chosen the classes that are presented in Table \ref{Classes}. 

\begin{table}[h]
\centering
\begin{tabular}{cl} 
Class & Age \\
\hline 1 & 0 - 15 \\
2 & 16 - 40 \\
3 & 41 - 70 \\
4 & 71+
\end{tabular}
\caption{Classes of patients of different age}
\label{Classes}
\end{table}

In this study, there are five groups of parameters:
\begin{enumerate}
\item parameters that are constant along the patients and not varied along the domain of computation,
\item parameters that are constant along the patients and varied along the domain of computation,
\item parameters that are varied along the patients and not varied along the domain of computation,
\item parameters that are varied along the patients and varied along the domain of computation,
\item parameters that are dependent on other parameters.
\end{enumerate}
To assess the uncertainty in the input data, we use a basic Monte Carlo method in which we sample input data from predefined statistical distributions. Regarding spatially heterogeneous parameters, we use sampling from a log-normal distribution. Each sample
is a one-dimensional realization, and is based on the heterogeneous sampling through a normalized truncated Karhunen-Lo{\'e}ve expansion of a zero-mean stochastic process, by
\begin{equation}\hat{u}(X)=\sum_{j=1}^{n} \hat{Z}_{j} \sqrt{\frac{2}{n}} \sin \left((2 j-1) \frac{\pi}{2\left|\Omega_{x, t}\right|} X\right),\end{equation}
where $\hat{Z}_{j} \sim \mathcal{N}(0,1)$, hence $\hat{Z}_{j}$ denotes a set of \emph{iid} stochastic variables that follow the standard normal distribution, $\left|\Omega_{x, t}\right|$ is the length of the domain of computation, and $-L \leq X \leq L$. From the stochastic variable $\hat{u}(X)$, we show the regeneration of, for example, $\hat{E}$ by
\begin{equation}\log (\hat{E}(X)) \sim \mu+\sigma \hat{u},\end{equation}
therewith
\begin{equation}\hat{E}(X)=\exp \left(\mathcal{M}_{E}+\mathcal{S}_{E} \hat{u}(X)\right).\end{equation}
Hence $\hat{E}(X)$ is a realization of a \emph{lognormal} distribution with mean $\mathcal{M}_E$ (expected value) and standard deviation $\mathcal{S}_E$. These values can be expressed by the arithmetic (sample) mean $\mu_E$ and arithmetic standard deviation $\sigma_E$ as follows
\begin{equation}
\mathcal{M}_{E}=\ln \left(\frac{\mu_{E}}{\sqrt{1+\frac{\sigma_{E}^{2}}{\mu_{E}^{2}}}}\right), \quad \mathcal{S}_{E}=\sqrt{\ln \left(1+\frac{\sigma_{E}^{2}}{\mu_{E}^{2}}\right)}.\end{equation}
In the same way, we can create heterogeneous, stochastic inputs for other parameters as well.

To test the model’s feasibility, we vary the parameter values based on the results found in literature on aging skin. We are interested in the differences in the intensity of contraction and the total strain energy between the distinct age classes. To quantify these differences, we test the null-hypothesis $H_0:\mu_A=\mu_B$ versus a two-sided alternative for classes $A$ and $B$ of patients using the following $t$-statistic:
$$t=\frac{\overline{X}_A-\overline{Y}_B}{s_p},\quad s_p = \sqrt{\frac{s_a^2 + s_b^2}{n_b}}$$
where $\overline{X}_A$ and $\overline{Y}_B$ are the mean values of the results in distinct age groups $A$ and $B$, $s_{p}$ is the estimated standard error of $\overline{X}_A-\overline{Y}_B$, $s_a^2$ and $s_b^2$ are the standard deviations in the age groups $A$ and $B$, and $n_b$ is the number of samples in the age groups. Here we assume that the number of samples in the age groups are equal. We reject the null-hypothesis if $|t|>t_{2(n_b-1)}(\alpha/2)$, with $\alpha=0.001$. 

\begin{figure}[h]
\centering
\includegraphics[width=\linewidth]{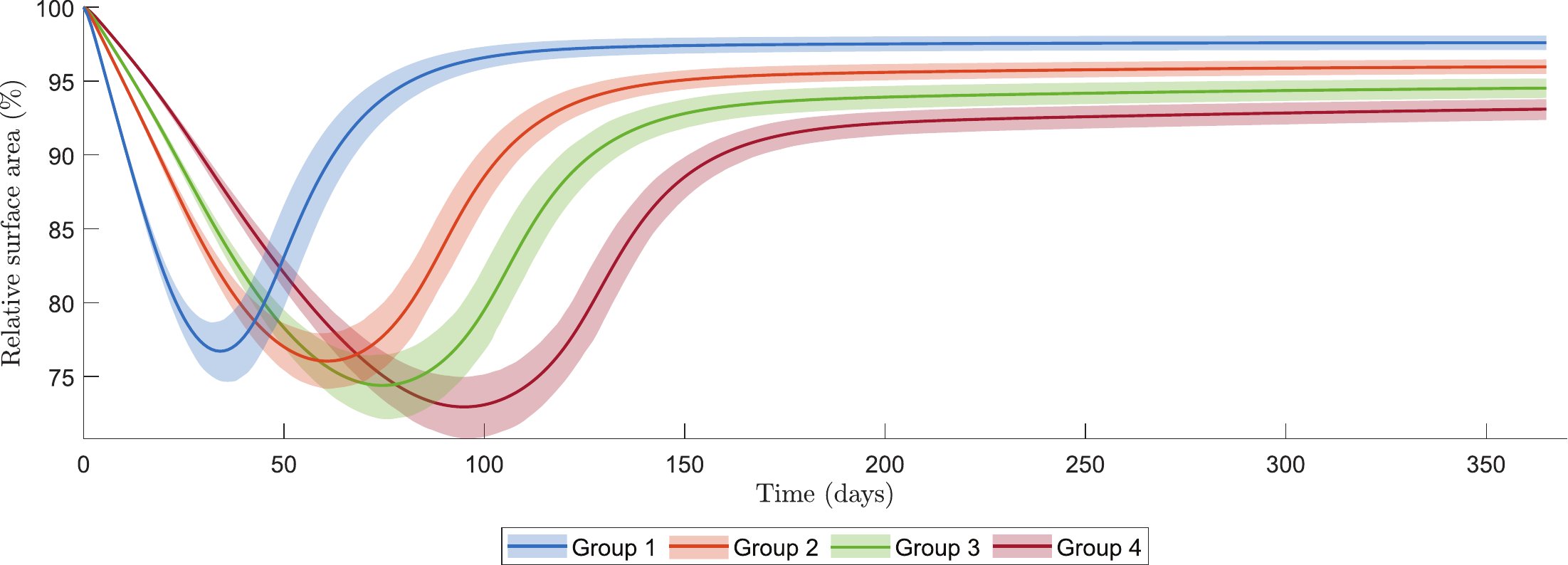}  
\caption{Confidence intervals for the contraction of burns in different age groups. The intervals show the mean values and the 95\% confidence values of the mean}
\label{fig:confidence}
\end{figure}

To reduce the computation time, we performed simulations on half a domain with $L = 10$ cm and $L^w = 3.6$ cm. 

We simulated $n_b=1950$ burns per age group. Hence, in total we simulated 7800 burns. We used parallel computing with 12 processors, three processors responsible per group, on a 64 bit Windows 10 Pro system with 16 GB RAM and 3.59 GHz AMD Ryzen 5 3600 6-Core Processor. The total computation time was 13.5 hours, hence the mean computation time per simulation takes less than half a minute. A major advantage of the one-dimensional implementation is its short computation time, which allows to do many of them within a reasonable time interval. 
For the test statistic, we used $t_{3898}(0.0005)=3.293$. The standard deviations $s_{i,m}^2$ for the age classes $i\in\{1,2,3,4\}$ for the minimum of the relative surface area (i.e., maximum contraction) are 
\begin{equation*}
\begin{aligned}
s_{1,m}^2 &= 1.2158, \quad s_{2,m}^2 = 1.1293, \quad s_{3,m}^2 = 1.2570,\\
s_{4,m}^2 &= 1.2815.
\end{aligned}
\end{equation*}
The standard deviations $s_{i,e}^2$ for the age classes $i\in\{1,2,3,4\}$ for the relative surface area on day 365 are 
\begin{equation*}
\begin{aligned}
s_{1,e}^2 &= 0.2941, \quad s_{2,e}^2 = 0.3105, \quad s_{3,e}^2 = 0.4003,\\
s_{4,e}^2 &= 0.4238.
\end{aligned}
\end{equation*}

\begin{figure}
\begin{subfigure}{.48\linewidth}
  \centering
  \includegraphics[width=\linewidth]{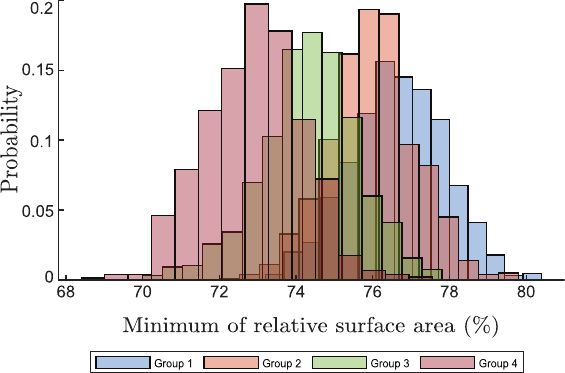}  
  \caption{}
  \label{fig:hist-rsaw-min}
\end{subfigure}\hfill
\begin{subfigure}{.48\linewidth}
  \centering
  \includegraphics[width=\linewidth]{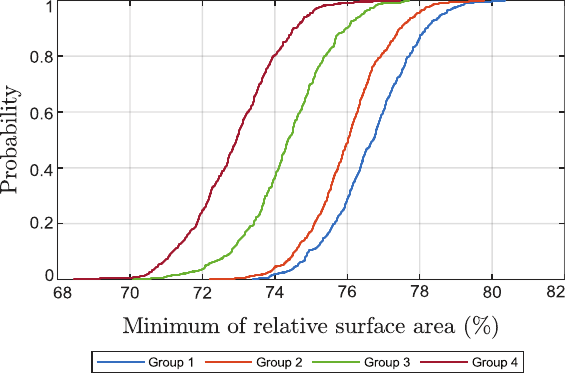}  
  \caption{}
  \label{fig:cdf-rsaw-min}
\end{subfigure}
\caption{Histograms of the minimum relative wound area (a) and its cumulative distribution function (b)}
\label{fig:rsaw-min}
\end{figure}

\begin{figure}[h]
\begin{subfigure}{.48\linewidth}
  \centering
  \includegraphics[width=\linewidth]{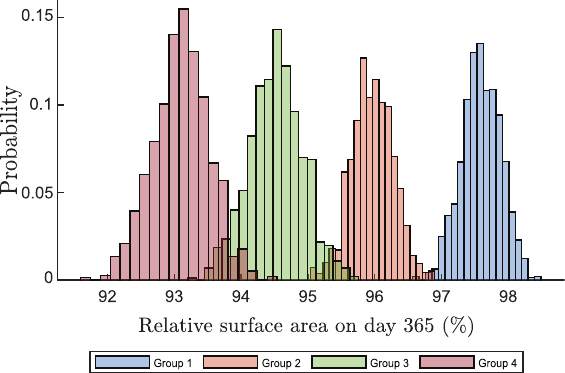}  
  \caption{}
  \label{fig:hist-rsaw-end}
\end{subfigure}\hfill
\begin{subfigure}{.48\linewidth}
  \centering
  \includegraphics[width=\linewidth]{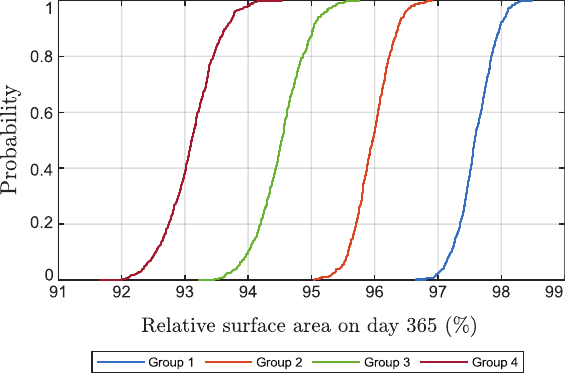}  
  \caption{}
  \label{fig:cdf-rsaw-end}
\end{subfigure}
\caption{Histograms of the relative surface area on day 365 (a) and its cumulative distribution function (b)}
\label{fig:rsaw-end}
\end{figure}

\begin{figure}[h]
\centering
\includegraphics[width=\linewidth]{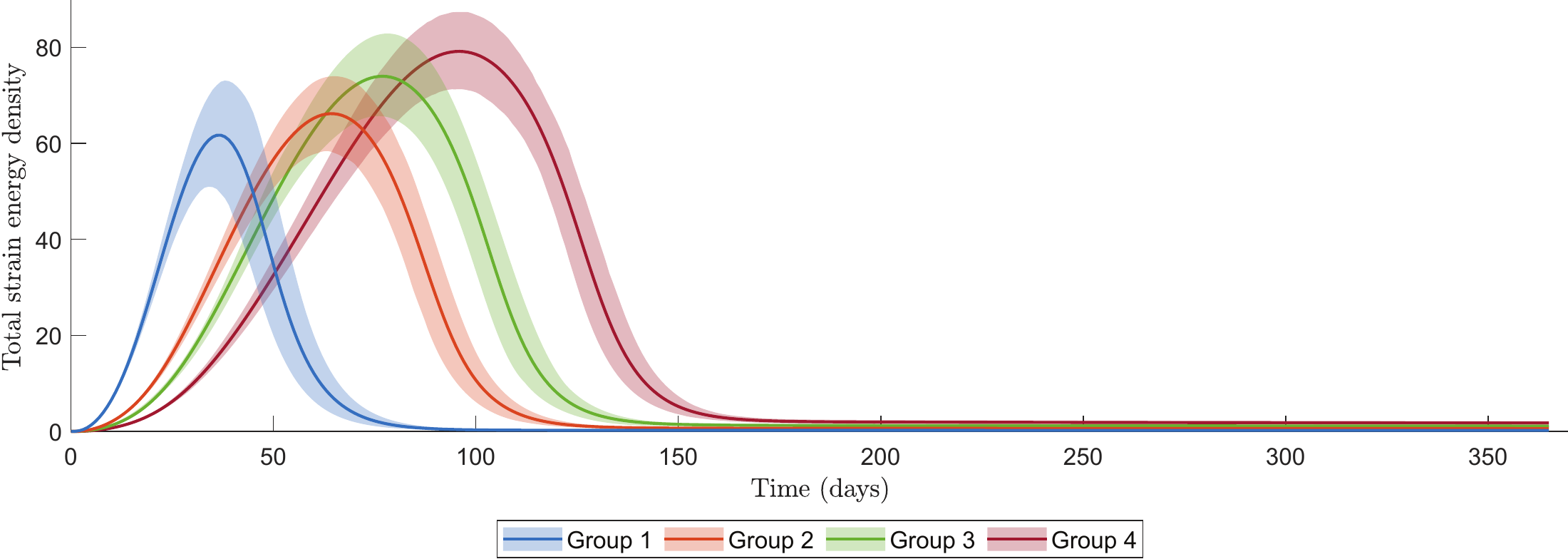}  
\caption{Confidence intervals for the strain energy in the healing of burns in different age groups. The intervals show the mean values and the 95\% confidence values of the mean}
\label{fig:confidence1}
\end{figure}

\begin{figure}
\begin{subfigure}{.48\linewidth}
  \centering
  \includegraphics[width=\linewidth]{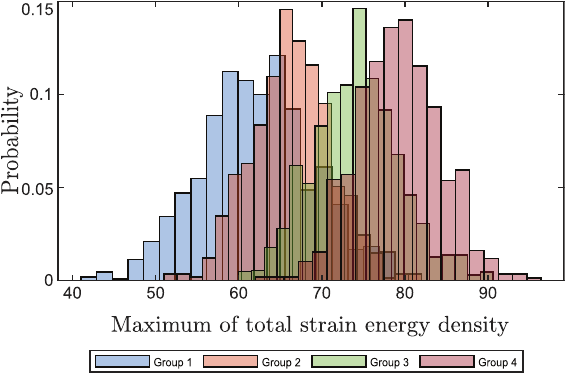}  
  \caption{}
  \label{fig:hist-strain-max}
\end{subfigure}\hfill
\begin{subfigure}{.48\linewidth}
  \centering
  \includegraphics[width=\linewidth]{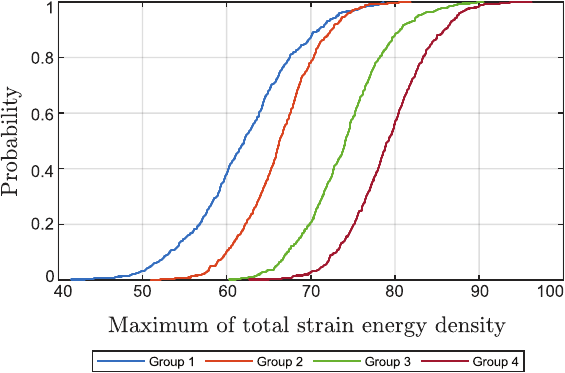}  
  \caption{}
  \label{fig:cdf-strain-max}
\end{subfigure}
\caption{Histograms of the maximum of the total strain energy density (a) and its cumulative distribution function (b)}
\label{fig:strain-max}
\end{figure}

To get some insight into the effects of the four classes of patients of different age shown in Table \ref{Classes}, we present confidence intervals for the contraction of burns in different age groups in Figure \ref{fig:confidence}. In Figures \ref{fig:rsaw-min} and \ref{fig:rsaw-end}, we show histograms of the minimum relative surface area and the relative surface area on day 365, respectively, together with the corresponding cumulative distribution function. Finally, we present confidence intervals for the total strain energy in the healing of burns in different age groups in Figure \ref{fig:confidence1}, and histograms of the maximum of the total strain energy density and its cumulative distribution function in Figure \ref{fig:strain-max}. The values for the parameters of the model, together with the standard deviation values for the variation over the domain of computation (if applicable), are shown in Table \ref{Parms} in the Appendix. 

Figure \ref{fig:confidence} shows four 95\% confidence intervals for the mean of the size of the scar, each confidence interval corresponding to a group of patients. The range of the contraction values comes from the variability of the parameters over the domain of computation. We can see that the maximum contraction value 
is about the same in the first two groups of ages, and from group 2, a higher age class gives a larger reduction of the size of the scar, and therewith a larger intensity of the contraction. Further, for higher ages, it takes more time to reach the maximum intensity of contraction. 

There seems to be more variability in the permanent deformation in the elderly patients. In the elderly patients, it takes longer before the wound healing cascade reaches equilibrium than in younger patients. Minima of the relative surface area were mostly reached on days 34, 61, 74, and 95, with values of 76.7, 76.0, 74.4, and 72.9\% for groups 1, 2, 3, and 4, respectively. Unfortunately, these results do not correspond fully to the observations in the clinic. Normally, contraction is of less order in elderly patients and, in general, retraction takes a longer period. We note that the longer retraction period is visible in the two-dimensional model results by Koppenol \cite{Koppenol2017a}, and that this ‘stretched retraction period’ is handled with the parameter $a_c^{II}$. Furthermore, in clinic one sees more contractures in younger patients in case the injury was in or near a joint. In contrast, in the clinic, the elderly seem to experience less discomfort because of contraction. One reason for this could be that the skin of the elderly is less tight than the skin of the young. Looser tissue can move more compared to skin that is already tight and is therefore less likely to cause movement restriction when it contracts.
The reason why our simulation results do not meet the clinical observations is because of the variety of factors that have not been modeled in our mathematical model yet.

Figure \ref{fig:rsaw-min} shows the results on the minimum relative surface area (i.e., maximum of contraction). Although we see that there is an overlap between all the groups in Fig. \ref{fig:hist-rsaw-min}, the maximal contraction is significantly different ($p<0.001$) between the groups (see Table \ref{tab:t-stat-min} in the Appendix). The differences in the minimum of the relative surface area between consecutive age groups are largest between ages 16-40 (group 2) and 41-70 (group 3), and smallest between ages 0-15 (group 1) and 16-40 (group 2). Given that in reality the evolution of the size of the scar is different from our results, we expect that the differences between group 1 and group 2 are larger in reality, possibly largest. 

The overlap between the age groups is also visible in the estimated cumulative distribution function plot in Fig. \ref{fig:cdf-rsaw-min}, where the functions of the first and second group, and the third and fourth group, almost intersect. From the cumulative distribution functions we can estimate the probabilities of reaching a certain amount of maximum contraction. For example, this figure suggests that with 70\% probability a patient from groups of patients 1, 2, 3, and 4, respectively, can reach 22.7\% contraction, 23.5\% contraction, 25.1\% contraction, and 26.4\% contraction. Using such functions in the future can help to predict the probability of developing a contracture, to intervene in time when possible. 

Figure \ref{fig:rsaw-end} shows the results on the relative surface area on day 365 (i.e., permanent contraction). Eventually, the scar maturates and because of the morphoelastic behavior of the skin, the size of the scar almost never reaches its initial configuration again. Here the initial configuration represents the initial size and geometry of the burn wound. Like for minimal contraction values, the overlap between consecutive groups for the intensity of contraction after one year in Fig. \ref{fig:hist-rsaw-end} are significant ($p<0.001$, see Table \ref{tab:t-stat-max} in the Appendix). The differences in the relative surface area after one year between consecutive groups are similar to the differences in the minimum of the relative surface area. The histograms show a correlation between aging and the spread, confirming the observation found in the confidence intervals, that there is larger variability in the intensity of contraction in elderly people than in children.

From the cumulative distribution functions in Fig. \ref{fig:cdf-rsaw-end}, we can estimate the probability of a certain contraction intensity. Here we see mean values for the final contraction about 97.4\%, 95.8\%, 94.5\%, and 93.3\% for groups 1, 2, 3, and 4, respectively. This final contraction intensity is an indicator for a possible contracture, in case the scar is in or near a joint. Given the location of the scar, one might estimate probabilities that this scar will develop a contracture of certain extent.

Figure \ref{fig:confidence1} shows four confidence intervals for the mean of the total strain energy density. The total strain energy density is a measure for the discomfort that a patient experience. Each confidence interval corresponds to a group of patients of different age. The range of the total strain energy density comes from the variability of the parameters over the domain of computation. We can see that the maximum of the total strain density is about the same in the first two groups, and from group 2, a higher age class gives a larger maximum of the total strain energy density. Further, for higher ages, it takes more time to reach the maximum total strain energy density. Note that all these results relate to the relative surface area densities shown in Figure \ref{fig:confidence}.

The total strain energy densities reach maxima on days 36, 64, 78, and 95, with values of 62, 66, 74 and 79 for groups 1, 2, 3, and 4, respectively. We note that the total strain energy density reaches maxima a few days later than the maximum contraction in almost all groups. Because the contraction data in Figure \ref{fig:confidence} does not relate to what doctors see in the clinic, we assume this is the same case for the total strain energy density. This means that in real-life children might experience more discomfort than the elderly, in contrast to what we see in Figure \ref{fig:confidence1}.

Figure \ref{fig:strain-max} shows the results on the maximum of the total strain energy density. Like in Figures \ref{fig:rsaw-min} and \ref{fig:rsaw-end}, although the histograms in Fig. \ref{fig:hist-strain-max} show an overlap between all the groups, the maximum of the total strain energy density is significantly different ($p<0.001$) between the groups (see Table \ref{tab:t-stat-max} in the Appendix).
This overlap is also visible in the cumulative distribution function plot in Fig. \ref{fig:cdf-strain-max}, where the functions of the first and second group intersect on the top, and the third and fourth group almost intersect. From the cumulative distribution functions we can estimate the probabilities of reaching a certain amount of maximum total strain. For example, this figure suggests that with 80\% probability a patient from groups of patients 1, 2, 3, and 4, respectively, can reach a maximum of 67.5, 70.4, 72, and 83.3 total strain. Since the figures show a strong correlation between the contraction and total strain, we conclude that targeting the contraction directly targets the total strain as well.

\section{Discussion}
\label{sec:7}
In this study, we worked with the morpho-elastic model for contraction and contractures in burn scars. This model was developed by Koppenol, and based on the principle of morpho-elasticity \cite{Hall}. We have provided a one-dimensional version of this model and focused on the parameter values, the sensitivity of the parameter values, and the feasibility of the model for patients of different ages. 

We comprehensively described the (ranges of the) values of the parameters. Most of the variety in the parameter values we have found in literature sources. We estimated some parameter values and adopted other parameter values from Koppenol. In case there were ranges of values found in literature, we chose upper or lower bounds, or a fixed value in between. We have used these values as mean values for the sensitivity analysis and for the parameters we have used for the second age group in our feasibility study. 

For our sensitivity analysis, we varied 30 parameters by $\pm 0,5,10,15,20,25$\%, and we varied the length of the initial wound. We showed results for the minimum of the relative surface area, the day on which wound healing reached this minimum, the relative surface area on day 365, the maximum value of the strain energy density, and the day on which wound healing reaches this maximum. The most sensitive parameter is the equilibrium collagen concentration present in the dermal layer. Other parameters that seem significantly sensitive are the apoptosis rate of fibroblasts, the apoptosis rate of myofibroblasts, the constant $R$ in equation \eqref{PDE_v} that influences the body force, and the secretion rate of signaling molecules. The parameters we found to be the least sensitive are the fluxes parameters, the crowding factor, the viscosity, the mass density of dermal tissue, and the initial collagen concentration. We note that we let the initial collagen concentration depend on the equilibrium collagen concentration, which can influence the sensitivity value of this parameter. 

Furthermore, we performed a feasibility study for the model to investigate the effect of aging on contraction, contractures and discomfort in burn wound healing.
We have chosen for four groups of patients in age classes: 0-15, 16-40, 41-70, 70+. 
We varied the parameters of the model according to observations from the literature, so that there was a variation between the groups of patients. 
We furthermore varied the parameters of the model using Karhunen-Lo{\'e}ve expansions, to model the heterogeneity of human skin, and in the framework of our Monte Carlo method we performed sampling from statistical distributions to assess the impact of uncertainty in the data on the behavior of contraction. 
The model is feasible for this approach, showing increased extent of contraction with age, a delayed maximum amount of contraction in elderly people (showing delayed healing), increased contracture in elderly people, and increased variety of contracture formation in elderly people compared to children. Next to these results, we see that the extent of discomfort is highly related to the contraction in wound healing. The figures show that there is a larger amount of discomfort in elderly patients and that the maximum discomfort is experienced significantly ($p<0.001$) earlier in younger children than in other age groups. 

This study shows that contraction increases with age and shows that there is a significant difference ($p<0.001$) in maximum amount of contraction between the different age groups. The least significant difference is found between ages 0-15 and 16-40, which can also be seen in Figure \ref{fig:confidence}. 
Further, the differences in the amount of contraction on day 365 in consecutive groups is least significant between 41-70 and 70+ years (see Table \ref{tab:t-stat-end} in the Appendix). The most significant difference in the amount of contraction on day 365 in consecutive groups is found to be between 0-15 and 16-40 years. Given the cumulative distribution functions, we can give probabilities of a certain amount of contraction in specific age-dependent groups of patients. We have seen that the differences in the maximum contraction and the contraction on day 365 are of a few percentages (less than 10) of order. For the maximum discomfort that a patient might feel, we have seen that there is significant difference ($p<0.001$) between all groups, of which in consecutive groups the difference between ages 16-40 and 41-70 is most significant. From the figures we can conclude that these patients experience the same amount of discomfort, although this happens much quicker in children. 

However, we note that we obtained the results with a mathematical model for which it is hard to find validation data. 
In the clinic, contraction seems to be of less order in elderly patients, and the retraction takes a longer period. 
The elderly have excess skin, meaning that they suffer less from contraction. 
There is less tension and less stretch in the skin of the elderly and hence we assume that there is less contraction. 
An explanation for our results might be that in fact such an extent of contraction is present in the elderly, however, this is not seen and noticed because of above reasons. 
On the other hand, the process of contraction is affected by growth and mobility.
These two kinds of forces yield a different process of contraction in agile, growing children compared to grown-ups. 
Therefore, we conclude that the model could be feasible for an age-related study for contraction in dermal wound healing and for that, the model needs further adjustments.

In this study, we defined the discomfort a patient may experience as the total strain energy. In equation \eqref{strainenergy}, we defined the total strain energy density for epsilon between $\Omega_t$ and $\Omega_t^{\text{eq}}$. However, it could well be that the definition where we take epsilon between $\Omega_{t=0}$ and $\Omega_t$ is more appropriate as a measure of the discomfort the patient is experiencing. We plan to consider this option in future studies, where we take into account two-dimensional effects.

Further research is needed to understand the differences of extent and timing of contraction in patients. To account for growth and movement forces, we plan on adding a factor for growth in equation \eqref{PDE_e}, and wound boundary forces for movement. To simulate for aging patients, we might put time-dependence on (some of) the parameters and simulate for a much longer time period. For example, the skin collagen content decreases at about 2\% per year \cite{Farage2015}, which has the most influence on the time healing period and amount of contraction. To distinguish between superficial and severe wounds, we can adapt real-life behavior in the initial conditions of the variables. In fact, in wounds where the environment still is vital, the release of cytokines is because cells lose their solid integrity. However, in severe wounds with a damaged environment, the release happens because of two reasons. Namely, the injury boiled signaling molecules, and the environment is dead. Here, circulation is absent and therefore little happens. Cytokines are present, but from vital edges where circulation still occurs, penetrating the tissue. The question is how deep, with an answer we can consider in the initial conditions.
To incorporate more geometrical matters, it is necessary to extend the model to a more-dimensional framework. 
Though more-dimensional frameworks allow to assess geometrical issues, these more-dimensional frameworks will require more simulation time and the use of a more advanced computer infrastructure if the objective is to carry out Monte Carlo simulations for the assessment of the likelihood that contraction of a particular intensity occurs. 
To reduce computational times, we want to model the boundaries of the wound as elastic springs. Further, it would be an improvement to code the finite element solution to the model in a high level programming language such as C++ and use an artificial intelligence framework such as neural networks. 
In finite element analysis, the necessary \emph{mesh generation} procedure is widely regarded as the weakest link in the chain of the analysis. To avoid the failure of a mesh to be analysis-suitable, we want to use isogeometric analysis (IGA) in the more-dimensional framework.
We plan to assess these issues in future work.

\section*{Acknowledgements}
The authors are grateful for the financial support by the Dutch Burns Foundation under Project 17.105.

\section*{Appendix}
\label{sec:8}
\begin{table}[h]
    \begin{subtable}[h]{0.40\textwidth}
	\begin{tabular}{ll} 
	Parameter & $\mu$  \\
	\hline\rule{0pt}{3ex}$\overline{M}$&0  \\
	$\overline{c}$	& 0 				 \\
	$\tilde{c}$ 		& $10^{-8}$ 	 \\
	$\delta_{M}$ 	& $6 \times 10^{-2}$ \\
	$R$ 			& 0.995 			\\
	$r_{F}^{\max }$ 	& 2 			\\
	$\eta^{I}$ 		& 2 			\\
	$\eta^{II}$ 		& 0.45 		\\
	$k_{\rho}^{\max }$ & 10 		
	\end{tabular}
       \caption{Constant along patients, constant along domain}
       \label{Parm1}
    \end{subtable}
    \hfill
    \begin{subtable}[h]{0.45\textwidth}
	\begin{tabular}{lll} 
	Parameter & $\mu$ & $\sigma$  \\
	\hline\rule{0pt}{3ex}$a_{c}^{I}$ 	& $10^{-8}$ 				& $3.45\times10^{-10}$ 	 \\
	$a_{c}^{I I}$ 	& $10^{-8}$ 				& $6.25\times10^{-10}$ 	 \\
	$a_{c}^{I V}$ 	& $10^{-9}$ 				& $10^{-10}$ 			 \\
	$\xi$ 			& $4.4 \times 10^{-2}$ 	& $1.1\times10^{-3}$ 		\\
	$\delta_{c}$ 		& $5 \times 10^{-4}$ 		& $9.8\times10^{-6}$ 		 \\
	$k_{c}$ 			& $3 \times 10^{-13}$ 		& $3.95\times10^{-15}$	\\
	$\rho_{t}$ 		& 1.09 					& $1.21\times10^{-1}$		\\
	$\chi_F$ 			& $2\times 10^{-3}$		& $2.22\times10^{-4}$		\\
			&							&
	\end{tabular}
        \caption{Constant along patients, varied along domain}
        \label{Parm2}
   \end{subtable}\\
   \begin{subtable}[h]{\textwidth}
   \begin{tabular}{lllll} 
	Parameter & $\mu^{1}$ & $\mu^{2}$ & $\mu^{3}$ & $\mu^{4}$ \\
	\hline\rule{0pt}{3ex}$\overline{N}$ 	& $1.5\times 10^{4}$ 	& $10^{4}$ 				& $9.5\times10^{3}$ 	& $9\times10^{3}$	 \\
	$\overline{\rho}$ 					& $1.25\times10^{-1}$	& $1.125\times10^{-1}$	& $1.05\times10^{-1}$	& $9.75\times10^{-2}$	
	\end{tabular}
	\caption{Varied along patients, constant along domain}
	\label{Parm3}   
   \end{subtable}\\
   \begin{subtable}[h]{\textwidth}
   \begin{tabular}{llllll} 
	Param. & $\mu^{1}$ & $\mu^{2}$ & $\mu^{3}$ & $\mu^{4}$ & $\sigma$  \\
	\hline\rule{0pt}{3ex}$D_{F}$	& $1.2\times10^{-6}$	& $10^{-6}$ 			& $8.5\times 10^{-7}$	& $7\times10^{-7}$	& $1.43\times10^{-7}$  \\
	$D_{c}$ 						& $3.25\times10^{-3}$ & $2.88\times10^{-3}$ & $2.55\times10^{-3}$	& $2.2\times10^{-3}$ 	& $3.2\times10^{-4}$	 \\
	$r_{F}$ 						& 1.222 				& $9.24\times10^{-1}$ & $8.16\times10^{-1}$	& $6.11\times10^{-1}$ & $7.11\times10^{-2}$	 \\
	$\kappa_{F}$ 				& $10^{-6}$ 			& $10^{-6}$ 			& $1.5\times10^{-6}$ 	& $2\times10^{-6}$ 	& $1.11\times10^{-7}$  \\
	$k_{F}$ 						& $1.14\times10^{7}$	& $1.08\times10^{7}$	& $1.02\times10^{7}$ 	& $9.63\times10^{6}$ 	& $5.68\times10^{5}$ 	\\
	$a_{c}^{III}$ 				& $2.05\times10^{8}$	& $2\times10^{8}$ 	& $1.95\times10^{8}$ 	& $1.9\times10^{8}$ 	& $4.35\times10^{6}$ 	 \\
	$\delta_N$					& $1.9\times10^{-2}$	& $2\times10^{-2}$	& $2.1\times10^{-2}$	& $2.2\times10^{-2}$	& $2.27\times10^{-4}$ \\
	$\delta_{\rho}$ 				& $6.11\times10^{-6}$ & $6\times10^{-6}$ 	& $5.89\times10^{-6}$ & $5.78\times10^{-6}$ & $1.09\times10^{-7}$  \\
	$\mu$ 						& $10^2$ 			& $10^2$ 			& $1.2\times10^2$ 	& $1.4\times10^2$ 	& $1.11\times10$ 	\\
	$E$ 							& $3.2\times10^2$	& $3.5\times10^2$ 	& $3.8\times10^2$	& $4.1\times10^2$	& $1.03\times10$	\\
	$\zeta$ 						& $3.8\times10^{2}$ 	& $4\times10^{2}$ 	& $4.2\times10^{2}$ 	& $4.4\times10^{2}$ 	& $1.82\times10$	
	\end{tabular}
	\caption{Varied along patients, varied along domain}
	\label{Parm4}
\end{subtable}
   \caption{Parameter values. In subtables (c) and (d), $\mu^i,i\in\{1,2,3,4\}$ denotes the mean value in class $i$ of patients of different age}
   \label{Parms}
\end{table}

\begin{table}[h]
    \begin{subtable}[h]{0.30\textwidth}
	\begin{tabular}{ll}
	Groups 	& $t$-value \\\hline
	1\&2 	& 17 \\
	1\&3	& 58 \\
	1\&4	& 94 \\
	2\&3	& 43 \\
	2\&4	& 80 \\
	3\&4	& 36
	\end{tabular}
	\caption{The values of the $t$-statistic comparing the minima of the relative surface area between age groups}
	\label{tab:t-stat-min}
    \end{subtable}
    \hfill
    \begin{subtable}[h]{0.30\textwidth}
	\begin{tabular}{ll}
	Groups 	& $t$-value \\\hline
	1\&2 	& 23 \\
	1\&3	& 63 \\
	1\&4	& 91 \\
	2\&3	& 49 \\
	2\&4	& 83 \\
	3\&4	& 32
	\end{tabular}
	\caption{The values of the $t$-statistic comparing the maxima of total strain between age groups}
	\label{tab:t-stat-max}
   \end{subtable}
   \hfill
   \begin{subtable}[h]{0.30\textwidth}
	\begin{tabular}{l|llllll}
	Groups 	&  $t$-value \\\hline
	1\&2 	& 167 \\
	1\&3	& 273 \\
	1\&4	& 385 \\
	2\&3	& 126 \\
	2\&4	& 241 \\
	3\&4	& 108 
	\end{tabular}
	\caption{The values of the $t$-statistic comparing the relative surface area on day 365 between age groups}
	\label{tab:t-stat-end} 
   \end{subtable}
   \caption{The values of the $t$-statistic comparing the minima of the relative surface area (a), the maxima of the total strain (b) and the relative surface area on day 365 (c)}
   \label{tab:t-stat}
\end{table}

\bibliographystyle{abbrv}
\bibliography{references}
\end{document}